\documentclass[11pt]{article}
\usepackage{}
\usepackage[total={6in, 9in}]{geometry}
\usepackage{amsfonts}
\usepackage{amsthm}
\usepackage{enumerate}
\usepackage{amssymb}
\usepackage{amsmath,amsfonts}
\usepackage{mathrsfs}
\usepackage{cases}
\usepackage{latexsym,bm}
\usepackage{indentfirst}
\usepackage{color}
\usepackage{ifpdf}
\usepackage{graphicx}
\usepackage{psfrag}
\usepackage{dsfont}
\usepackage[
pdfauthor={CESS},
pdftitle={Spanning trails},
pdfstartview=XYZ,
bookmarks=true,
colorlinks=true,
linkcolor=blue,
urlcolor=blue,
citecolor=blue,
bookmarks=true,
linktocpage=true,
hyperindex=true
]{hyperref}

\newtheorem{THM}{\textbf{Theorem}}

\newtheorem{LEM}{\textbf{Lemma}}[section]
\newtheorem{CLA}{\textbf{Claim}}

\newcommand{\pf}{\noindent\textbf{Proof}.\quad}

\newenvironment{claimproof}[1]{\par\noindent\underline{Proof:}\space#1}{\hfill $\blacksquare$}

\newcommand{\iC}{\overset{\leftharpoonup }{C}}

\newcommand{\oC}{\overset{\rightharpoonup }{C}}

\begin{document}
\title{Hamiltonian cycles  in tough $(P_2\cup P_3)$-free graphs}
\author{%
	Songling Shan \vspace{0.5cm}\\ 
	 Illinois State University, Normal, IL 61790\\
		\texttt{sshan12@ilstu.edu}\\
}

\date{8 January 2019}
\maketitle
\maketitle

\emph{\textbf{Abstract}.}
Let $t>0$ be a real number and $G$ be a graph. We say $G$
is $t$-tough if for every cutset $S$ of $G$, the ratio of 
$|S|$ to the number of components of $G-S$ is at least $t$.  Determining toughness is an NP-hard problem for arbitrary graphs. The Toughness Conjecture of Chv\'atal, stating that there exists a constant $t_0$ such that every $t_0$-tough graph with at least three vertices  is hamiltonian, is still open in general.  
A graph is called $(P_2\cup P_3)$-free if it does not contain any induced subgraph isomorphic to $P_2\cup P_3$, the union of two vertex-disjoint paths of order 2 and 3, respectively.  In this paper, we show that every 15-tough $(P_2\cup P_3)$-free graph with at least three vertices is hamiltonian. 

%
%


\vspace{2mm}

\section{Introduction}

Graphs considered in this paper are simple, undirected,  and finite.
Let $G$ be a graph.
Denote by $V(G)$ and  $E(G)$ the vertex set and edge set of $G$,
respectively. For $v\in V(G)$,  $N_G(v)$ denotes the set of neighbors
of $v$ in $G$. 
For $S\subseteq V(G)$  and $x\in V(G)$,
define $deg_G(x, S)=|N_G(x)\cap S|$.  
If $H\subseteq G$, we simply write $deg_G(x, H)$ for  $deg_G(x, V(H))$. 
We skip the subscript $G$ if the graph in consideration is clear from the context. 
Let $S\subseteq V(G)$.  Then the subgraph induced on $V(G)-S$ is denoted by
$G-S$. For notational simplicity, we write $G-x$ for $G-\{x\}$.
If $uv\in E(G)$ is an edge, we write $u\sim v$.
Let $V_1,
V_2\subseteq V(G)$ be two disjoint vertex sets. Then $E_G(V_1,V_2)$ is the set
of edges of $G$ with one end in $V_1$ and the other end in $V_2$.

The number of components of $G$ is denoted by $c(G)$. Let $t\ge 0$ be a
real number. The graph is said to be {\it $t$-tough\/} if $|S|\ge t\cdot
c(G-S)$ for each $S\subseteq V(G)$ with $c(G-S)\ge 2$. The {\it
toughness $\tau(G)$\/} is the largest real number $t$ for which $G$ is
$t$-tough, or is  $\infty$ if $G$ is complete. This concept, a
measure of graph connectivity and ``resilience'' under removal of
vertices, was introduced by Chv\'atal~\cite{chvatal-tough-c} in 1973.
 It is  easy to see that  if $G$ has a hamiltonian cycle
then $G$ is 1-tough. Conversely,
 Chv\'atal~\cite{chvatal-tough-c}
 conjectured that
there exists a constant $t_0$ such that every
$t_0$-tough graph is hamiltonian (Chv\'atal's toughness conjecture).
 Bauer, Broersma and Veldman~\cite{Tough-counterE} have constructed
$t$-tough graphs that are not hamiltonian for all $t < \frac{9}{4}$, so
$t_0$ must be at least $\frac{9}{4}$.

There are a number of papers on
 Chv\'atal's toughness conjecture,
and it  has
been verified when restricted to a number of graph
classes~\cite{Bauer2006},
including planar graphs, claw-free graphs, co-comparability graphs, and
chordal graphs. A
graph $G$ is called {\it $2K_2$-free\/} if it does not contain two independent
edges as an induced subgraph.
Recently, Broersma, Patel and Pyatkin~\cite{2k2-tough} proved that
every 25-tough $2K_2$-free graph on at least three vertices is
hamiltonian, and the author of this paper improved the required 
toughness in this result from 25 to 3~\cite{1706.09029}.

Let $P_\ell$ denote a path on $\ell$-vertices. 
A graph is $(P_2\cup P_3)$-free if it does not contain any induced 
copy of $P_2\cup P_3$, the disjoint union of  $P_2$ and $P_3$.    
In this paper, we confirm Chv\'atal's toughness conjecture
for the class of $(P_2\cup P_3)$-free graphs, a superclass 
of $2K_2$-free graphs.

\begin{THM}\label{main}
Let $G$ be a $15$-tough $(P_2\cup P_3)$-free graph with at least three vertices. Then $G$
is hamiltonian.
\end{THM}

In~\cite{MR1392734} it was shown that every 3/2-tough split graph on at least three vertices
is hamiltonian.  And the authors constructed
a sequence $\{G_n\}_{n=1}^{\infty}$ of split graphs with no 2-factor and $\tau(G_n)\rightarrow 3/2$.
So $3/2$ is the best possible toughness for split graphs to be hamiltonian.
Since split graphs are $(P_2\cup P_3)$-free,
we cannot decrease the bound in Theorem~\ref{main} below 3/2. Although it is certain that 15-tough is not optimal, 
we are not sure about the ``best possible''  toughness for giving a hamiltonian cycle in a  $(P_2\cup P_3)$-free graph. 

The class of $2K_2$-free graphs is well studied, for instance, see~\cite{
	2k2-tough, CHUNG1990129, MR845138, 2k21,  MR2279069,  MR1172684}.
It is a superclass of {\it split\/} graphs,
where the vertices can be partitioned into a clique and an independent set.
One can also easily check that every {\it cochordal\/} graph (i.e., a graph that is the complement of a
chordal graph) is $2K_2$-free and so the class of $2K_2$-free graphs is
at least as rich as the class of chordal graphs.  
By the definition, the class of $(P_2\cup P_3)$-free graphs 
is a superclass of $2K_2$-free graphs but with  much more complicated structures than graphs that are $2K_2$-free. 
The proof techniques used in~\cite{2k2-tough} and~\cite{1706.09029}
for showing that certain tough $2K_2$-free graphs are hamiltonian seem to be not 
applicable for $(P_2\cup P_3)$-free graphs. The proof approach used in this paper for showing Theorem~\ref{main}
is new and more general and  reveals some 
structural properties of $(P_2\cup P_3)$-free graphs. 



\section{Proof of Theorem~\ref{main}}

We start this section with some definitions.
Let $G$ be a graph and $S\subseteq V(G)$ a cutset of $G$, and let $D$ be a component of $G-S$. For a vertex $x\in S$, we say that {\it$x$ is adjacent to $D$\/} if $x$ is adjacent in $G$ to a vertex of  $D$. 
We call $D$ a {\it clique component\/} of $G-S$
if $V(D)$ is a clique in $G$.  We call $D$ a {\it trivial component\/} of $G-S$ if $D$ has only one vertex. 

A {\it star-matching\/} is a set of vertex-disjoint copies of stars.  The vertices of degree at least 2
in a star-matching are called the centers of the star-matching.  In particular, if all 
the stars in a {\it star-matching\/} are isomorphic to $K_{1,t}$, where $t\ge 1$ is an integer, we call the star-matching a 
{\it $K_{1,t}$-matching\/}.   For a star-matching $M$, we denote by  $V(M)$
the set of vertices covered by $M$. 

 Let $C$ be an oriented cycle. For $x\in V(C)$,
denote the successor of $x$ on $C$ by $x^+$ and the predecessor of $x$ on $C$ by $x^-$.
For $u,v\in V(C)$, $u\oC v$ denotes the portion of $C$
starting at $u$, following $C$ in the orientation,  and ending at $v$.
Likewise, $u\iC v$ is the opposite portion of $C$ with endpoints as $u$
and $v$.   We  assume all cycles
in consideration afterwards are oriented. 
A path $P$ connecting two vertices $u$ and $v$ is called 
a {\it $(u,v)$-path}, and we write $uPv$ or $vPu$ in specifying the two endvertices of 
$P$. Let $uPv$ and $xQy$ be two paths. If $vx$ is an edge, 
we write $uPvxQy$ as
the concatenation of $P$ and $Q$ through the edge $vx$.

\begin{LEM}[\cite{MR2816613}, Theorem 2.10]\label{matching}
	Let $G$ be a bipartite graph with partite sets $X$ and $Y$, and let $f$ be a function from $X$ to the set of
	positive integers. If for every $S\subseteq X$,  it holds that 
	$|N_G(S)|\ge \sum\limits_{v\in S}f(v)$ ,  then $G$ has a subgraph $H$ such that $X\subseteq V(H)$,
	$d_H(u)=f(u)$ for every $u\in X$, and $d_H(v)=1$ for every $v\in Y\cap V(H)$. 
	\end{LEM}

\begin{LEM}[Bauer et al.~\cite{MR1336668}]\label{degree-tough}
Let $t>0$ be real and $G$  be a $t$-tough $n$-vertex graph ($n\ge 3$)  with  $\delta(G) > \frac{n}{t+1}-1$.  Then 
 $G$ is hamiltonian. 
\end{LEM}

Lemmas~\ref{clique-components} and~\ref{clique-components2} below  are consequences of $(P_2\cup P_3)$-freeness. 

\begin{LEM}\label{clique-components}
Let $G$ be a  $(P_2\cup P_3)$-free graph and 
$S\subseteq V(G)$ a cutset of $G$. If  $G-S$  has a component that is not a clique component, 
then all other components of $G-S$ are trivial.  Consequently, if $G-S$
has at least two nontrivial components, then all components of $G-S$ are clique components. 

\end{LEM}	

\begin{LEM}\label{clique-components2}
	Let $G$ be a  $(P_2\cup P_3)$-free graph and 
	$S\subseteq V(G)$ a cutset of $G$, and let $x\in S$. 
	Suppose that $x$ is adjacent to exactly one component $D$ of 
	$G-S$, and $G-S$ has a nontrivial component to which $x$
	is not adjacent, then $x$ is adjacent in $G$ to all vertices of $D$. 
	\end{LEM}	

\begin{LEM}\label{clique-connection}
Let $G$ be a connected $(P_2\cup P_3)$-free graph and 
$S\subseteq V(G)$ a cutset of $G$ such that each vertex in $S$
is adjacent to at least two components of $G-S$. 
Then each of the following statement holds.
\begin{enumerate}[(i)]
	\item For every nontrivial clique component $D\subseteq G-S$ and for every 
	vertex $x\in S$, $x$ is adjacent to $D$.  
	\item For every nontrivial clique component $D\subseteq G-S$ and for every 
	vertex $x\in S$,  if $x$ is adjacent to at least three components, then $x$ is adjacent in $G$ to at least $|V(D)|-1$ vertices of $D$. 
	\item Let $D_1$ and $D_2$ be two  nontrivial clique components of $G-S$.  Then for  every 
	vertex $x\in S$,  either $x$ is adjacent in $G$ to at least $|V(D_i)|-1$ vertices of each $D_i$, or $x$
	is adjacent in $G$ to all vertices of one of $D_i$,  $i=1,2$. 
\end{enumerate}
\end{LEM}
	
	\pf  For (i), let $w_1$ and $w_2$ be two neighbors of $x$ in $G$
	respectively from two distinct components of $G-S$.
	Then $w_1xw_2$ is an induced $P_3$. Now for every edge $uv\in E(D)$, 
	we must have that $\{w_1, w_2\}\cap \{u,v\}\ne \emptyset$ or $x$ is adjacent in $G$ to $u$ or $v$,
	by the $(P_2\cup P_3)$-freeness. 
	Therefore, $x$ is adjacent to $D$. 
	For (ii), let $x\in S$
	and $D$ be a nontrivial clique component of $G-S$. Since $x$ is adjacent to at least three components, there exists $u,w$, respectively 
	from  two components of $G-S$ that are distinct from $D$ such that $x\sim u$
	and $x\sim w$ in $G$. Thus, $uxw$ is an induced $P_3$ in $G$.
	Furthermore, since $u,w\in V(G)-S-V(D)$, $E_G(\{u,w\}, V(D))=\emptyset$. 
	Thus, by   the $(P_2\cup P_3)$-freeness assumption, for every edge in $D$, $x$
	is adjacent to at least one endvertex of that edge. This, together with the fact that $D$
	is a clique, we know that $x$ is adjacent at least $|V(D)|-1$ vertices in $D$. 
	For (iii), assume to  the contrary that the statement does not hold. 
	By symmetry, we assume that there exists $uv\in E(D_1)$ such that $x\not\sim u,v$ in $G$, and there exists $w\in V(D_2)$ such that $x\not\sim w$ in $G$.  Let $y\in V(D_2)\cap N_G(x)$ that exists by Lemma~\ref{clique-connection} (i). Then $uv \cup xyw$
	is an induced $P_2\cup P_3$, giving a contradiction. 
	\qed

\begin{LEM}\label{cycle-extendabilit y}
	Let $t\ge 1$  and $G$ be an $n$-vertex  $t$-tough graph, 
	and let $C$ be a non-hamiltonian cycle of $G$. 
If $x\in V(G)-V(C)$ satisfies  that $deg(x, C)>\frac{n}{t+1}$,  then $G$ has a cycle $C'$ such that $V(C')=V(C)\cup \{x\}$. 	
\end{LEM}

\pf It is clear that if $x$ is adjacent to two consecutive vertices $u,w$ on $C$, then  
$$C'=(C-\{uw\})\cup \{ux,xw\}$$ 
is a cycle with the desired property. So we assume that for any $u,w\in N_G(x)\cap V(C)$, $uw\not\in E(C)$. Let $W=\{u^+\,|\, u\in N_G(x)\cap V(C)\}$ be the set of the successors of the neighbors of 
$x$ on $C$.  Because there is a one-to-one correspondence between $W$ and $N_G(x)\cap V(C)$,
by the assumption that $deg(x, C)>\frac{n}{t+1}$, we know that 
\begin{equation}\label{size of W}
|W|>\frac{n}{t+1}.
\end{equation}
If there exist $u^+, w^+\in W$ with $u, w\in N_G(x)\cap V(C)$ such that $u^+ \sim w^+$ in $G$, 
then $$C'=u^+\oC wxu\iC w^+u^+$$
is a desired cycle. Therefore, we assume that $W$ is an independent set in $G$. 
Let $S=V(G)-W$.  Then $c(G-S)=|W|$, as $W$ is an independent set in $G$. However,  by~\eqref{size of W}, 
$$
\frac{|S|}{c(G-S)}=\frac{|S|}{|W|}<\frac{\frac{tn}{t+1}}{\frac{n}{t+1}}=t,
$$
showing a contradiction to the toughness of $G$. 
\qed 

\begin{LEM}\label{cycle-extendability2}
	Let  $G$ be an $n$-vertex  $15$-tough $(P_2\cup P_3)$-free graph,  
	and let $C$ be a non-hamiltonian cycle of $G$. 
	Let $P \subseteq G-V(C)$ be an $(x,z)$-path.  If both $x$ and $z$ are adjacent to 
	more than $\frac{4.5n}{16}$ vertices in $V(C)$, then $G$ has a cycle $C'$ such that $V(C')=V(C)\cup V(P)$. 	
\end{LEM}

\pf It is clear that if $x$ is adjacent to a vertex $u$ on $C$ and $z$ is adjacent  to a vertex $w$ on $C$
such that $uw\in E(C)$, then  
$$C'=(C-\{uw\})\cup \{ux,zw\}\cup P$$ 
is a cycle with the desired property. So we assume that 
\begin{equation}\label{non-consecutive}
\mbox{for any $u\in N_G(x)\cap V(C)$
	and any $w\in N_G(z)\cap V(C)$, $uw\not\in E(C)$.}
\end{equation}
 Let 
\begin{eqnarray*}
W_x&=&\{u^+\,|\, u\in N_G(x)\cap V(C)\},\\
W_z&=&\{u^+\,|\, u\in N_G(z)\cap V(C)\}.
\end{eqnarray*}
 Clearly, 
 \begin{equation}\label{Wxzsize}
|W_x|=|N_G(x)\cap V(C)|>\frac{4.5n}{16}, \quad \mbox{and} \quad |W_z|=|N_G(z)\cap V(C)|> \frac{4.5n}{16}. 
 \end{equation}
If there exist $u^+\in W_x$ and $w^+\in W_z$ with $u\in N_G(x)\cap V(C)$  
and $w\in N_G(z)\cap V(C)$ such that $u^+ \sim w^+$ in $G$, 
then $$C'=u^+\oC wzPxu\iC w^+u^+$$
is a desired cycle. Therefore, we assume that 
\begin{equation}\label{noedge}
E_G(W_x,W_z)=\emptyset.
\end{equation}

We further claim that 
\begin{equation}\label{noncons2}
\mbox{no two vertices in $N_G(x)\cap V(C)$  or $N_G(z)\cap V(C) $ are consecutive on $C$.}
\end{equation}
By symmetry, we only show that no two vertices in $N_G(x)\cap V(C)$ are consecutive on $C$.

Assume to the contrary that there exists a path $v_1v_2\cdots v_\ell\subseteq C$  with $\ell\ge 2$ such that for each $i$ with $1\le i\le \ell$, $v_i\in N_G(x)\cap V(C)$, $v_1^-\not\in N_G(x)\cap V(C)$, and $v_\ell^+\not\in N_G(x)\cap V(C)$. 
Note that such vertices $v_1$ and   $v_\ell$ exist by the assumption in~\eqref{non-consecutive}. 
We assume that 
there exist $w_1,w_2\in W_z$ such that  $w_1\sim w_2$ in $G$. 
For otherwise, letting $S=V(G)-W_z$ gives that $c(G-S)=|W_z|$, and consequently,  $|S|/c(G-S)=|S|/|W_z|<15$ by~\eqref{Wxzsize}. 

Then  $xv_\ell v_\ell^+$ is an induced $P_3$
in $G$.   Consider the edge $w_1w_2$. 
By the assumption in~\eqref{non-consecutive}, $x\not\sim w_1,w_2$ in $G$ (otherwise, $w_1^-w_1\in E(C)$ or $w_2^-w_2\in E(C)$ with $w_1^-, w_2^- \in N_G(z)\cap V(C)$),
and by the assumption in~\eqref{noedge}, $v_\ell^+\not\sim w_1,w_2$
in $G$. Thus, $v_\ell\sim w_1$ or $v_\ell \sim w_2$ in $G$ by the $(P_2\cup P_3)$-freeness assumption. 
However, $v_\ell=v_{\ell-1}^+\in W_x$, showing a contradiction to~\eqref{noedge}. 

Therefore, by~\eqref{noncons2}, 
\begin{equation}\label{disjoint}
(N_G(x)\cap V(C))\cap W_x=\emptyset, \quad \mbox{and}\quad (N_G(z)\cap V(C))\cap W_z=\emptyset. 
\end{equation}

Let 
$$
W_{xz}=W_x\cap W_z.
$$
By the assumption in~\eqref{noedge}, $W_{xz}$ is an independent set in $G$. By the 
toughness of $G$, we know that $|W_{xz}|\le \frac{n}{16}$. 
Therefore,  $|N_G(x)\cap N_G(z)\cap V(C)|\le \frac{n}{16}$. 
These, together with~\eqref{Wxzsize} and~\eqref{disjoint}, imply that 
\begin{eqnarray*}
n&\ge &|(N_G(x)\cap V(C))\cup (N_G(z)\cap V(C))\cup W_x\cup W_z|\\&>& \frac{9n}{16}+\frac{9n}{16}-|W_{xz}|-|N_G(x)\cap N_G(z)\cap V(C)|\\
&\ge& \frac{16n}{16}=n,
\end{eqnarray*}
showing a contradiction. 
\qed 

%

\begin{LEM}\label{cycle-extendability3}
	Let $G$ be an $n$-vertex  $15$-tough $(P_2\cup P_3)$-free graph with $n\ge 31$,
	and let $S\subseteq V(G)$ be a cutset of $G$ with $|S|\le \frac{3n}{4}$.
	Assume that  $G-S$ has at least two nontrivial clique components, and  that for every edge $uv\in E(G)$, $d(u)+d(v)\ge |S|$. 
	 Then 
	 $G$ has a hamiltonian cycle.  
\end{LEM}

\pf  By Lemma~\ref{clique-components}, every component of $G-S$ is a clique component. 
If there exists $x\in S$ such that $x$ is adjacent to exactly one component, say $D$
of $G-S$, then we move $x$ from $S$ into $D$. 
By Lemma~\ref{clique-components2},  every  component of 
$G-(S-\{x\})$ is still  a clique component.   We move out all such vertex 
$x$ from $S$ iteratively and denote the remaining vertices in $S$ 
by $S_1$.  Note that $S_1\ne \emptyset$, since $G$ is a connected graph and 
$S$ is a cutset of $G$.  Also, $c(G-S)=c(G-S_1)$ and $G-S_1$ has at least two nontrivial components.
By Lemma~\ref{clique-components}, every component of $G-S_1$ is a clique component. 
Let 
\begin{eqnarray*}
S_0&=&\{x\in S_1\,|\, \mbox{$x$ is not adjacent to any component of $G-S_1$}\},\\
S_2&=&\{x\in S_1\,|\, \mbox{$x$ is  adjacent to at least two  components of  $G-S$}\}. 
\end{eqnarray*}
Note that $S_2=S_1-S_0$.

Since $G-S_1$ has a nontrivial component that has no edge going to $S_0$, 
the $(P_2\cup P_3)$-freeness of $G$ implies that $G[S_0]$ consists of vertex-disjoint complete subgraphs of $G$. 
Thus $S_2$
is a cutset of $G$ with components consisting those from $G-S_1$ and $G[S_0]$. 
Also,  all 
components of $G-S_2$ are clique components in which at least two of them are nontrivial.  
By the toughness of $G$, $|S_2|\ge 15c(G-S_2)$.

We will construct a hamiltonian cycle in $G$ through two steps: (1) combing spanning cycles from every clique component of $G-S_2$
that has at least three vertices into a single 
cycle $C$, and (2) insterting remaining vertices  in $V(G)-V(C)$ into $C$
to obtain a hamiltonian cycle of $G$.  

Suppose that $G-S_2$ has exactly $h$ clique components $D_1,\cdots, D_h$ 
with $|V(D_1)|\ge |V(D_2)|\ge \cdots \ge |V(D_h)|\ge 1$, and that 
the first $t$ ($0\le t\le h$) of them are components that contain at least three vertices. 
	Since $G-S_2$ has at least two nontrivial components, both $D_1$ and $D_2$ are nontrivial. 

\begin{CLA}\label{D1}
The component $D_1$ contains at least 5 vertices. 
\end{CLA}
\begin{claimproof}
	Since $|S|\le \frac{3n}{4}$, $n\ge \frac{4|S|}{3}$. Also,  $c(G-S)\le \frac{|S|}{15}$ by $\tau(G)\ge 15$. Therefore, 
	a largest component of $G-S$ contains at least
	\begin{eqnarray*}
\frac{n-|S|}{c(G-S)} &\ge & \frac{\frac{4|S|}{3}-|S|}{\frac{|S|}{15}}\ge 5
	\end{eqnarray*}
 vertices. 
\end{claimproof}

Let
\begin{eqnarray*}
	Q_1&=&\{x\in S_2\,|\, \mbox{$x$ is adjacent to a component distinct from $D_1$ and $D_2$}\},\\
	Q_2&=&\{x\in S_2\,|\, \mbox{$x$ is adjacent to less than $\frac{|V(D_1)|-1}{2}$ vertices of $D_1$}\},\\
	Q_3&=&\{x\in S_2\,|\, \mbox{$x$ is adjacent to less than $\frac{|V(D_2)|-1}{2}$ vertices of $D_2$}\}.
\end{eqnarray*}
By Lemma~\ref{clique-connection} (i) and the definition of $Q_1$, we know that if $Q_1\ne \emptyset$, then every vertex in 
$Q_1$ is adjacent to at least three components of $G-S_2$. By Lemma~\ref{clique-connection} (ii), 
we get the following claim. 
\begin{CLA}\label{Q1}
Suppose that $Q_1\ne \emptyset$. Then  for every $x\in Q_1$ and for every nontrivial component $D$ of $G-S$,   $x$ is adjacent to at least $|V(D)|-1$ vertices of  $D$.
\end{CLA}

\begin{CLA}\label{Q2}
	Suppose that $Q_2\ne \emptyset$.   Then for  every $x\in Q_2$, $x$ is adjacent to all vertices of  $D_2$ and $Q_2$ is a clique in $G$.
\end{CLA}
\begin{claimproof}
Note that both $D_1$ and $D_2$ are nontrivial components of $G-S_2$.  
Since $|V(D_1)|\ge 5$ by Claim~\ref{D1},  $x$ is not adjacent to at least three 
vertices of $D_1$ by the definition of $Q_2$. Therefore, 
$x$ is adjacent to all vertices  of  $D_2$ by  Lemma~\ref{clique-connection} (iii). 
For the second part, 
	suppose to the contrary that there exist $x,y\in Q_2$ such that $x\not\sim y$ in $G$.
	Let $w\in V(D_2)$. Then $w\sim x$ and $w\sim y$ in $G$ by the first part of this claim. Thus, 
	we find an induced $P_3=xwy$.
	Since $E_G(\{w\}, V(D_1))=\emptyset$, 
	the  $(P_2\cup P_3)$-freeness  implies that for every edge in $D_1$,
	at least one of $x$ and $y$ is adjacent to at least one endpoint of the edge. Since $D_1$ is complete,  by Pigeonhole Principle, one of $x$ and 
	$y$ is adjacent to at least $\frac{|V(D_1)|-1}{2}$ vertices of $D_1$.  This  gives a contradiction to the assumption that $x,y \in Q_2$. 
\end{claimproof}

Similarly,  we have the following result. 
\begin{CLA}\label{Q3}
Suppose that $Q_3\ne \emptyset$.   Then for every $x\in Q_3$, $x$ is adjacent to all vertices  of  $D_1$ and $Q_3$ is a clique in $G$.
\end{CLA}
By Claims~\ref{Q1} to \ref{Q3}, we have that 
\begin{equation}\label{Qi-relation}
 Q_i\cap Q_j=\emptyset, i\ne j,  i,j=1,2,3.
\end{equation}
Define 
\begin{eqnarray*}
	W&=&\bigcup_{\max\{t+1, 3\}\le i\le h} V(D_i).
\end{eqnarray*}
If $W\ne \emptyset$, we claim that there is a $K_{1,2}$-matching $M$ between $W$ and $S_2$ 
such that every vertex in $W$ is the center of a $K_{1,2}$-star. 
For otherwise, by  Theorem~\ref{matching},  there exists $W_1\subseteq W$ such that $2|W_1|>|N_G(W_1)\cap S_2|$.
By the definition of $W$,  we then see that $G-(N_G(W_1)\cap S_2)$ has at least $|W_1|/2+1\ge 2$
components (vertices in $S_2-(N_G(W_1)\cap S_2)$ form at least one component of $G-(N_G(W_1)\cap S_2)$ that is disjoint from those containing vertices from $W_1$), implying that 
$$\frac{|N_G(W_1)\cap S_2|}{c(G-(N_G(W_1)\cap S_2))}<4<15.$$
This gives a contradiction to the toughness.
\begin{equation}\label{defM}
\mbox{Let $M$ be a $K_{1,2}$-matching between $W$ and $S_2$.}
\end{equation}

\begin{CLA}\label{Mvertex}
	If $M\ne \emptyset$, then for every $x\in V(M)\cap S_2$, $x\in Q_1$.  Consequently,  for every nontrivial component $D$ of $G-S$,   $x$ is adjacent to at least $|V(D)|-1$ vertices of  $D$. 
\end{CLA}

\begin{claimproof} 
If $G-S_2$ has at least three nontrivial components, then
 every vertex of $S_2$ is adjacent to all those nontrivial components by 
Lemma~\ref{clique-connection} (i). 
Therefore, $S_2=Q_1$  by the definition of $Q_1$.
In particularly,  $x\in Q_1$ for  $x\in V(M)\cap S_2$. 
Hence, we assume that $G-S_2$ has exactly two nontrivial components, which are $D_1$ and $D_2$. 
This assumption implies that $|V(D_3)|\le 1$.
Consequently, $|V(D_3)|=1$ since $M\ne \emptyset$.  
Then for every $x\in V(M)\cap S_2$,
$x$ is adjacent to both $D_1$ and $D_2$ by Lemma~\ref{clique-connection} (i). Also $x$ is adjacent to 
a trivial component of $G-S_2$. Thus 
 $x\in Q_1$.  The second part of Claim~\ref{Mvertex} is a consequence of Claim~\ref{Q1}.
\end{claimproof} 
	
\begin{CLA}\label{HcycleC}
	There is a cycle $C$ in $G-V(M)$ with at least $\frac{3n}{20}$ vertices such that $C$ contains 
	all vertices from every $D_i$, $i=1,2,\cdots, t$, and $Q_2\cup Q_3\subseteq V(C)$. 
\end{CLA}

\begin{claimproof} 
		Suppose first that $G-S_2$ has at least three nontrivial components. 
		Then by Lemma~\ref{clique-connection} (i),   every vertex of $S_2$ is adjacent to 
		all those nontrivial components of $G-S_2$. Consequently, $S_2=Q_1$ and $Q_2=Q_3=\emptyset$.  
		 Therefore, for every 
	$x\in S_2$ and every $D_i$, $x$ is adjacent to at least $|V(D_i)|-1$ vertices of $D_i$ by Claim~\ref{Q1}. 

Let $x_1,\cdots, x_t$ be $t$ distinct vertices in $S_2-V(M)$. (By the toughness of $G$, $|S_2|\ge 15c(G-S_2)$. Since $|V(M)\cap S_2|\le 4c(G-S_2)$,  we have enough vertices in $S_2-V(M)$ to pick.)
	Let $C_i$ be a hamiltonian cycle of $D_i$, and let $u_i, v_i\in V(C_i)$ with $u_iv_i\in E(C_i)$ 
such that for $i=1,2,\cdots, t-1$, $x_i\sim v_i, u_{i+1}$,
and $x_t\sim u_1, v_t$ in $G$. Then $$C=u_1\oC_1v_1x_1u_2\oC_2v_2\cdots u_{t-1}\oC_{t-1}v_{t-1}x_{t-1}u_t\oC_tv_tx_tu_1$$
is a cycle that contains all vertices from each $D_i$ and  the vertices 
$x_1,\cdots, x_t$ from $S_2-V(M)$.	 Also $Q_2\cup Q_3\subseteq V(C)$ trivially as $Q_2=Q_3=\emptyset$. 

So we assume that 
$G-S_2$ has exactly two nontrivial clique components,  which  are $D_1$ and $D_2$, call this {\bf assumption ($\ast$)}.  
Let 
$$G_1=G[V(D_1)\cup Q_3] \quad \mbox{and}\quad  G_2=G[V(D_2)\cup Q_2].$$
By Claims~\ref{Q2} and~\ref{Q3}, we know that both $G_1$
and $G_2$ are complete subgraphs of $G$.   

If $t=1$ and $Q_2=\emptyset$, we let $C$ be a hamiltonian cycle of $G_1$.  Clearly, $ Q_2\cup Q_3\subseteq V(C)$.

Thus, we assume that $t\ge 2$ or $Q_2\ne \emptyset$. 
Since $D_1$ has at least three vertices by Claim~\ref{D1}, $G_1$ contains at least three vertices.  
Note that $|V(D_2)|\ge 2$
by the assumption that $G-S_2$ has at least two nontrivial components and $D_2$
is one of them.  
 Thus,  $G_2$ contains at least three vertices either by $t\ge 2$ or $Q_2\ne \emptyset$.  

If there are two disjoint edges between $G_1$
and $G_2$, then $G[V(G_1)\cup V(G_2)]$ has  a hamiltonian cycle $C$. Thus, 
we may assume, without loss of generality, that there is either no edge between $G_1$ and 
$G_2$ or all edges between $G_1$ and $G_2$ are incident to only a single vertex, say in $G_1$.

If $c(G-S_2)=2$, then  $M=\emptyset$ by the definitions of $W$ and $M$.  
Since $G[V(G_1)\cup V(G_2)]$ 
has a cutvertex or is disconnected,  the toughness of $G$ implies that $|S_2-Q_2-Q_3|\ge 29$. 
In addition, there are vertex-disjoint paths $P_1$ and $P_2$  connecting $G_1$
and $G_2$ in $G$ such that each $P_i$ only has exactly one of its endvertices in $G_1$ 
and $G_2$.  Let $V(P_i)\cap V(G_1)=\{x_i\}$ and  $V(P_i)\cap V(G_2)=\{y_i\}$, $i=1,2$. 
Let $C_1$ be a hamiltonian cycle in $G_1$ such that $x_1x_2\in E(C_1)$, and 
$C_2$ be a hamiltonian cycle in $G_2$ such that $y_1y_2\in E(C_2)$. Then 
$$C=x_1P_1y_1\oC_2y_2P_2x_2\iC_1x_1$$
is a cycle that contains all vertices in clique components of $G-S_2$ that contain at least three vertices and the vertices  from $P_1$ and $P_2$. Also  $Q_2\cup Q_3\subseteq V(C)$ by the construction of $C$.

So we assume that $c(G-S_2)\ge 3$.  Since  $D_1$ and $D_2$
are  the only nontrivial components of $G-S_2$ by assumption ($\ast$), the assumption that $c(G-S_2)\ge 3$ 
implies that  $|Q_1|\ge 15(c(G-S_2)-2)$. Thus $|Q_1-V(M)|\ge 15(c(G-S_2)-2)-2(c(G-S_2)-2)\ge 13$, since each trivial component of $G-S$ uses exactly two vertices from  $S_2\cap V(M)$.  Hence, we can find two vertices 
$x,y\in Q_1-V(M)= S_2-Q_2-Q_3-V(M)$ such that both $x$ and $y$
are adjacent to at least $|V(D_1)|-1$ vertices of  $D_1$, and at least   $|V(D_2)|-1$ vertices of $D_2$ by Claim~\ref{Mvertex}. 
We claim that $x$ is adjacent to at least two vertices  
of  $G_2$. This is clear if  $x$ is adjacent to at least two vertices of $D_2$. 
So we 
assume that $x$ is adjacent to only one 
vertex of $D_2$.  Let $w$  be the neighbor of $x$ from $D_2$, and let $w_1\in V(D_2)-\{w\}$, $w_2\in V(G_2)-V(D_2)$. By this choice, $x\not\sim w_1, w_2$ in $G$.    Note that $w_1$ is not adjacent to any vertex of $D_1$, and $w_2$
is adjacent to less than $\frac{|V(D_1)|-1}{2}$ vertices of  $D_1$. Therefore, we 
can find a vertex $w^*\in V(D_1)$ such that $w_1,w_2\not\sim w^*$ in $G$ and $x\sim w^*$ in $G$.   By the 
choice of $x$, there is a vertex $w'\in V(G)-S-V(D_1)-V(D_2)$ such that $x\sim w'$ in $G$. 
However, $w_1w_2\cup w^*xw'$ is an induced $P_2\cup P_3$. This gives a contradiction. 
Since $D_1$ has at least 5 vertices, both $x$ and $y$ have at least four neighbors in $D_1$. 
Thus we can select distinct vertices $x_1,y_1\in V(G_1)$
 and  $x_2,y_2\in V(G_2)$ such that 
 $x\sim x_1,x_2$ and  $y\sim y_1,y_2$ in $G$. 

Let $C_1$ be a hamiltonian cycle of $G_1$ such that  $x_1y_1\in E(C_1)$, and let $C_2$ be a hamiltonian cycle of $G_2$ such that  $x_2y_2\in E(C_2)$. Then 
$$
C=x_1xx_2\oC_2y_2yy_1\iC_1x_1
$$
is a cycle that contains all vertices in clique components of $G-S_2$ that contain at least three vertices and the vertices $x$ and $y$.  Furthermore, $Q_2\cup Q_3\subseteq V(C)$. 

Since for each $i$, $1\le i\le t$, $V(D_i)\subseteq V(C)$ and $\bigcup_{1\le i\le t} V(D_i)\subseteq V(G)-S-W$, we have that 
\begin{eqnarray*}
|V(C)|&\ge& n-|S|-|W| \ge n-|S|-2c(G-S)\\
&\ge & n-|S|-\frac{2|S|}{15}\ge n-\frac{17}{15}\cdot \frac{3n}{4}=\frac{3n}{20}. 
\end{eqnarray*}
\end{claimproof}	

\begin{CLA}\label{morethan_n_over_16}
	Let $C$ be the cycle defined in Claim~\ref{HcycleC}. 
For any $x\in S_2-V(C)$, $x$ has more than $\frac{n}{16}$ neighbors on  $C$. 
\end{CLA}

\begin{claimproof}
	Note that every vertex in 
	$S_2$ is adjacent to at least two components of $G-S_2$. 
	If $G-S_2$ has at least three nontrivial clique components, then 
	Lemma~\ref{clique-connection} (ii) implies that for 
	every  $x\in S_2$, and for every nontrivial clique component $D$
	of $G-S_2$, $x$ is adjacent to at least $|V(D)|-1$
	vertices of $D$.  
	By the toughness of $G$ and the assumption that $|S|\le \frac{3n}{4}$, $x$ has at least 
$$
	n-|S|-\frac{2|S|}{15}\ge \frac{3n}{20} 
$$
	neighbors in the union of the nontrivial clique components of $G-S_2$
	that contain at least three vertices.  Therefore,  $x$ has more than $\frac{n}{16}$ neighbors on  $C$. 
	
	So we assume that $G-S_2$ has exactly two nontrivial clique components. Since $S_2-V(C)\subseteq  S_2-Q_2-Q_3$ (recall that $Q_2\cup Q_3\subseteq V(C)$), we know that 
	$x$ is adjacent to at least $\frac{|V(D_1)|-1}{2}$ vertices of $D_1$, and is adjacent to at least $\frac{|V(D_2)|-1}{2}$ vertices of  $D_2$.  Note that 
	$$
	|V(D_1)|+|V(D_2)|\ge n-|S|-\frac{|S|}{15}\ge \frac{n}{5}. 
	$$
	Since $C$ contains all vertices from  $D_1\cup D_2$, we conclude that $x$
	is adjacent to at least $\frac{n}{10}-1>\frac{n}{16}$  (by $n\ge 31$) neighbors on  $C$. 
\end{claimproof}

By Claim~\ref{morethan_n_over_16}, and by 
applying  Lemma~\ref{cycle-extendabilit y}  for $C$ and vertices in $S_2-V(C)-V(M)$  iteratively, 
we  get a longer cycle $C'$ such that $V(C')=V(C)\cup (S_2-V(C)-V(M))$. 
Note also that 
$$
S_2-V(C')=V(M)\cap S_2, \quad \text{and} \quad  V(G)-S_2-V(C') =V(M)\cap (V(G)-S_2).
$$

Recall that for every $x\in S_2-V(C')=S_2\cap V(M)$,  $x$ is adjacent to 
at least $|V(D_i)|-1$ vertices in each $D_i$, $i=1,2,\cdots, t$ by Claim~\ref{Mvertex}. Therefore, if $|S|\le \frac{7n}{12}$,  then $x$ has in $G$ at least 
$$
\frac{5n}{12}-\frac{2|S|}{15}\ge\frac{5n}{12}-\frac{14n}{15\cdot12} >\frac{n}{3}>\frac{4.5n}{16}$$
neighbors  on $C'$.  Then applying Lemma~\ref{cycle-extendability2}  for  $C'$ and   every path in $M$ iteratively, 
we obtain a hamiltonian cycle in $G$. 
Hence  we
\begin{equation}\label{Slowersize}
\text{Assume that $|S|> \frac{7n}{12}$}.
\end{equation}

\begin{CLA}\label{halfless}
	For any two $K_{1,2}$-stars $x_1u_1y_1, x_2u_2y_2\in M
	$, if $u_1u_2$
	is a 2-vertex component of $G-S_2$ and $|S|>\frac{7n}{12}$, then at least one of $u_1$ and $u_2$ has more than $\frac{n}{16}$
	neighbors on  $C'$. 
\end{CLA}

\begin{claimproof}
	For otherwise,  since $u_i$ is adjacent to 
	exactly one vertex in $V(M)\cap (V(G)-S_2)$, and $|V(M)\cap S_2|\le 4|V(M)\cap (V(G)-S_2)|\le \frac{4|S|}{15}$, 
	$$
	d_G(u_1)+d_G(u_2)\le 2\left(\frac{n}{16}+1+|V(M)\cap S_2|\right)\le 2\left(\frac{n}{16}+1+\frac{4|S|}{15}\right)<|S|, 
	$$
	show a contradiction to the assumption that for every edge $uv\in E(G)$, $d_G(u)+d_G(v)\ge |S|$.  
\end{claimproof}

 Let 
\begin{eqnarray*}
M_1=\{uwv\in M\,|\, deg_G(w, C')>\frac{n}{16} \}, &\quad M_2=M-M_1. 
\end{eqnarray*}
Take  $uwv\in M_1$, note that  $u, v\in S_2$ and  $w\in  V(G)-S_2$. 
 By the definition of $M_1$, $deg(w, C')>\frac{n}{16}$. 
 By Claim~\ref{morethan_n_over_16},  $deg(u, C')>\frac{n}{16}$  and $deg(v, C')>\frac{n}{16}$. 
Now 
applying  Lemma~\ref{cycle-extendabilit y} for $C'$ and every path in $M_1$ iteratively, 
we get a longer cycle $C^*$ such that $V(C^*)=V(C')\cup V(M_1)$.

By the toughness of $G$, $G-S_2$
has at most $\frac{|S|}{15}$ components in total. Particularly,  $G-S_2$
has at most $\frac{|S|}{15}$ components that have at most two vertices in total.
By Claim~\ref{halfless}, we know that  for every 2-vertex component $uv$ of $G-S_2$, 
at least one of $u$ or $v$ has more than $\frac{n}{16}$ neighbors on $C'$. 
Therefore, at least one of the two $K_{1,2}$-stars centered, respectively, at $u$ and $v$
is contained in $M_1$.  In other words, there is at most one $K_{1,2}$-star
from $M_2$ that centers at a vertex from a same component of $G-S_2$.  Therefore, 
$$|V(M_2)|\le \frac{|S|}{15}+ \frac{2|S|}{15}=\frac{|S|}{5}.$$

By the definition of $M_2$ and 
by the assumption that for any $uv\in E(G)$, $d_G(u)+d_G(v)\ge |S|$, we know that  for any 
path $xwy\in M_2$, where $x,y\in V(G)-S_2$ 
and $w\in S_2$, we have that  $d_G(x)+d_G(w)\ge |S|$. 
Therefore,   the number of neighbors that $x$
has in $G$  on $C^*$ is at least 
\begin{eqnarray*}
&&|S|-deg_G(x,  G-V(C^*))-d_G(w)\\
&\ge & |S|-deg_G(x, V(M_2))-\big(deg_G(w, C^*)+deg_G(w, S_2\cap V(M_2))\big)\\
&\ge & |S|-\frac{|S|}{5}-\left(\frac{n}{16}+\frac{2|S|}{15}\right)=\frac{2|S|}{3}-\frac{n}{16}\\
&>&\frac{2\cdot 7n}{3 \cdot 12}-\frac{n}{16}=\frac{7 n}{24}>\frac{4.5n}{16}.
\end{eqnarray*}
 Similarly, the vertex $y$ 
has in $G$ at least  $\frac{4.5n}{16}$
neighbors on  $C^*$. 
Now applying 
Lemma~\ref{cycle-extendability2} for  $C^*$ and  every path in $M_2$ iteratively gives a hamiltonian cycle in $G$.  
\qed

\proof[Proof of Theorem~\ref{main}]
Since $G$ is 15-tough, it is 30-connected,
and consequently, $\delta(G)\ge 30$. 
By Lemma~\ref{degree-tough},
we may assume that 
\begin{equation}\label{lowern}
n\ge (\delta(G)+1)\cdot (\tau(G)+1)\ge 31\cdot 16, \quad \mbox{and}\quad \delta(G)\le \frac{n}{16}-1. 
\end{equation}
We consider two case to finish the proof. 

{\bf \noindent Case 1: For every edge $e=uv\in E(G)$, $d_G(u)+d_G(v)>\frac{3n}{4}$}.

Denote by  
\begin{equation}\label{V1-def}
V_1=\{v\in V(G)\,|\, d_G(v)\le \frac{3n}{8}\}.
\end{equation}
By the assumption of Case 1, we know that $V_1$
is an independent set in $G$. Therefore, 
\begin{equation}\label{V1-size}
|V_1|\le \frac{n}{16}.
\end{equation}
by $\tau(G)\ge 15$. 

Since $G$ is 15-tough, Lemma~\ref{matching} implies that 
$G$ has a $K_{1,2}$-matching  $M$ with all vertices in $V_1$ 
as the centers of the $K_{1,2}$-matching. 
Let $V_2$
be the set of the vertices contained in $M$. 
By~\eqref{V1-size}, we have that 
\begin{equation}\label{V2-size}
|V_2|\le \frac{3n}{16}.
\end{equation}

Denote by $G_1=G-V_2$. Then by the definitions  of $V_1, V_2$
and \eqref{V2-size}, we 
get that 
\begin{eqnarray}\label{G1-deg}
\delta(G_1)&>&\frac{3n}{8}-|V_2|\ge \frac{3n}{16}, \label{G1-deg} \\
deg_G(x, G_1)&>& \frac{3n}{8}-|V_2|\ge \frac{3n}{16}, \quad \mbox{for any $x\in V_2-V_1$}. \label{V2-V1-deg}
\end{eqnarray}


We first assume that $G_1$ has a hamiltonian cycle $C$.  
For every copy of $K_{1,2}$, say $xyz\in M$,  by~\eqref{V2-V1-deg},

\begin{eqnarray}
deg_G(x, G_1)>\frac{3n}{16}>\frac{n}{16},\nonumber\\
\label{xzdegree}\\
deg_G(z, G_1)>\frac{3n}{16}>\frac{n}{16}.\nonumber
\end{eqnarray}

Let 
\begin{eqnarray*}
	M_1=\{uwv\in M\,|\, deg_G(w, C)>\frac{n}{16} \}, &\quad M_2=M-M_1. 
\end{eqnarray*}
  By~\eqref{xzdegree}, 
applying  Lemma~\ref{cycle-extendabilit y} with respect to $C$ and every vertex in $M_1$ iteratively, 
we get a longer cycle $C^*$ such that $V(C^*)= V(C)\cup V(M_1)$.

By the definition of $M_2$ and 
by the assumption that for any $uv\in E(G)$, $d_G(u)+d_G(v)> \frac{3n}{4}$, we know that  for any 
path $xwy \in M_2$, where $x,y\in V_2-V_1$ 
and $w\in V_1$, we have that  $d_G(x)+d_G(w)> \frac{3n}{4}$. 
Therefore,   the number of neighbors that $x$
has in $G$  on $C^*$ is at least 
\begin{eqnarray*}
	&& \frac{3n}{4}-deg_G(x,  G-V(C^*))-d_G(w)\\
	&\ge &  \frac{3n}{4}-deg_G(x, V(M_2))-\big(deg_G(w, C^*)+deg_G(w, V_2)\big)\\
	&\ge &  \frac{3n}{4}-|V_2|-\left(\frac{n}{16}+|V_2-V_1|\right)\\
	&\ge &  \frac{3n}{4}-\frac{3n}{16}-\frac{n}{16}-\frac{2n}{16}\\
	&=&\frac{6n}{16}>\frac{4.5n}{16}.
\end{eqnarray*}
Similarly, the vertex $y$ 
has in $G$ at least  $\frac{4.5n}{16}$
neighbors on  $C^*$. 
Now applying 
Lemma~\ref{cycle-extendability2} for  $C^*$ and  every path in $M_2$ iteratively gives a hamiltonian cycle in $G$.

Hence we assume that $G$ does not have a hamiltonian cycle.  By Lemma~\ref{degree-tough} and~\eqref{G1-deg}, we know that $\tau(G_1)<7$. Therefore, there exists $S_1\subseteq V(G_1)$
such that 
\begin{equation}\label{toughG1}
|S_1|/c(G_1-S_1)<7.
\end{equation}
If $|S_1|\ge \frac{3n}{16}$,  then we get  that $c(G_1-S_1)=c(G-(S_1\cup V_2))$, and thus by~\eqref{V2-size}, 
$$
\frac{|S_1\cup V_2|}{c(G-(S_1\cup V_2))}<14, 
$$
showing a contradiction to $\tau(G)\ge 15$. So we assume that $|S_1|<\frac{3n}{16}$, thus 
$|S_1|\le \lfloor \frac{3n}{16}\rfloor$. 
As $\delta(G_1)\ge \lfloor \frac{3n}{16}\rfloor+1$ by~\eqref{G1-deg}, we know that each component of 
$G_1$ contains at least 
$$\delta(G_1)-|S_1|\ge \lfloor \frac{3n}{16}\rfloor+1-\lfloor \frac{3n}{16}\rfloor+1=2
$$
vertices.  By Lemma~\ref{clique-components}, we know that every component of 
$G_1-S_1$ is a clique component.   Let $S=S_1\cup V_2$. 
We then see that all components of $G-S$ are nontrivial.  Also, $|S|<\frac{6n}{16}<\frac{3n}{4}$ 
since  $|S_1|<\frac{3n}{16}$ and   $|V_2|\le \frac{3n}{16}$ by~\eqref{V2-size}.  
Furthermore, by the assumption of Case 1, for every edge $uv\in E(G)$, $d_G(u)+d_G(v)> \frac{3n}{4}>|S|$. 
Now we can apply  Lemma~\ref{cycle-extendability3} on $G$ and $S$ to find a 
 hamiltonian cycle in $G$.

{\bf \noindent Case 2: There exists an edge $e=uv\in E(G)$ such that  $d_G(u)+d_G(v)\le \frac{3n}{4}$}.

Let 
$$
S=N_G(u)\cup N_G(v)-\{u,v\},
$$
such that $d_G(u)+d_G(v)$ is smallest among all the degree sums of two adjacent vertices in $G$. 

By the assumption of this case and the choice of $S$, we know that 
\begin{equation}\label{S-size}
|S|\le \frac{3n}{4}-2, \quad \mbox{and} \quad \mbox{for any $u'v'\in E(G)$, $d(u')+d(v')\ge |S|$.}
\end{equation}
By the definition of $S$, $c(G-S)\ge 2$ and $uv$ is one of the components of $G-S$. 
Since $\tau(G)\ge 15$, and $|V(G)-S-\{u,v\}|\ge n-|S|\ge  \frac{|S|}{3}=\frac{5|S|}{15}$,
$G-S-\{u,v\}$ has a component with at least 5 vertices.  This, together with the fact that 
$uv$ is one of the components of $G-S$, Lemma~\ref{clique-components} implies that 
every component of $G-S$ is a clique component, and 
$G-S$ has at least two nontrivial components.   Again Lemma~\ref{cycle-extendability3} implies that $G$
has a hamiltonian cycle. 
\qed

\bibliographystyle{plain}
\bibliography{SSL-BIB}

\end{document}